
\documentstyle{amsppt}
\NoRunningHeads
\NoBlackBoxes
\magnification = 1200

\define\a{\alpha}
\predefine\barunder{\b}
\redefine\b{\beta}
\define\bpage{\bigpagebreak}

\define\cd{\cdot}
\predefine\dotunder{\d}
\redefine\d{\delta}
\define\e{\epsilon}

\define\f{\frac}
\define\g{\gamma}

\define\lb\{{\left\{}

\define\lm{\limits}

\define\mpage{\medpagebreak}

\define\oper{\operatorname}

\define\rb\}{\right\}}
\define\s{\sigma}

\define\sub{\subheading}

\define\th{\theta}

\define\({\left(}
\define\){\right)}
\define\[{\left[}
\define\]{\right]}
\define\<{\left<}
\define\>{\right>}
\def\slantline#1#2#3#4#5{\hbox to 0pt{
}}



\def\BR{{\bold R}}

\def\BZ{{\bold Z}}


\define\sos#1#2{(S_{#1}, #2)}
\define\vecs#1#2{#1_1,\dots,#1_#2}
\define\tec{Teichm\"uller\ }
\define\tecs{Teichm\"uller's\ }
\define\Ree{\oper{Re}}

\define\im{\oper{Im}}
\define\qd{\oper{QD}}
\define\ext{\oper{ext}}

\topmatter
\title Teichm\"uller Space is Not Gromov Hyperbolic
\endtitle
\thanks *Partially supported by the NSF grant DMS 9201321\endgraf
**Partially supported by the NSF grants DMS 9300001 and DMS 9022140 (MSRI);
Alfred
P\.
Sloan Research Fellow\endthanks
\endtopmatter

\document

\centerline{Howard A. Masur*}

\centerline{Department of Mathematics}

\centerline{University of Illinois at Chicago}

\centerline{Chicago, IL 60680}

\bpage
\centerline{Michael Wolf**}

\centerline{Department of Mathematics}

\centerline{Rice University}

\centerline{Houston, TX 77251}

\bpage
\centerline{Address until 8/1/94: Mathematical Sciences Research
Institute}

\hskip143pt 1000 Centennial Drive

\hskip143pt  Berkeley, CA 94720

\vskip.7cm
The \tec space of surfaces of genus $g>1$ with the \tec metric is
not nonpositively curved, in the sense that there are distinct
geodesic rays from a point that always remain within a bounded
distance of each other (\cite{Ma1}.) Despite this phenomenon,
\tec
space and its quotient, Moduli space, share many properties with
spaces of negative curvature: for instance, most converging
geodesic rays are asymptotic \cite{Ma2}, and the geodesic flow on
the moduli space is ergodic \cite{Ma3}.

One can ask whether these properties can be explained by \tec
space
having non-positive curvature in a sense weaker than that of
Busemann
used in \cite{Ma1}, which declared a space $X$ to be negatively
curved if the endpoints of two segments from $p\in X$ are spread
more than twice as far as the midpoints.

In this study of hyperbolic groups, Gromov (\cite{Gr}, see also
\cite{GdlH}) introduced a notion of negative curvature, now
called
Gromov hyperbolicity, that still captured many of the qualitative
aspects of Riemannian negative sectional curvature, but was less
restrictive than that of Busemann. Specifically, Gromov declared
a
space $X$ to be hyperbolic if there existed a number $M$ so that
for any $p\in X$ and any triangle in $X$ with vertex at $p$, the
leg of the triangle opposite $p$ would be within an
$M$-neighborhood of the legs of the triangle emanating from $p$.
Thus, for instance, the flat Euclidean strip
$\{(x,y)\in\BR^2\mid0<x<1\}$ would be Gromov hyperbolic but not
Buseman negatively curved; moreover, the fact that there are
pairs
of rays emanating from $p\in T_g$ which do not diverge does not,
in
itself, preclude \tec space with the \tec metric from being
Gromov
hyperbolic.

Nevertheless, the goal of this paper (Theorem~3.1) is to show
that
\tec space is not Gromov hyperbolic. This, of course, also
immediately implies that any quasi-isometric metric to the \tec
metric is also not Gromov hyperbolic, so any Gromov hyperbolic
metric on the \tec space is quite different from the \tec metric.

In this connection, one needs to observe that the isometry group
of
the \tec metric is the mapping class group (\cite{Roy}),
which contains large rank abelian subgroups, and so is not a
Gromov
hyperbolic group (with the word metric).
This in itself
does not seem to imply immediately that \tec space is not Gromov
hyperbolic.
For example there are Kleinian groups
 with rank $2$ abelian subgroups acting on hyperbolic $3$ space,
a Gromov
 hyperbolic space. It
 does suggest that good
  candidates for triangles to
contradict
Gromov's condition might be constructed with vertices at images
of a
single point $p$ under high iterates of commuting isometries.

In fact, this is the approach we take, showing (Theorem~3.1) that
with respect to the Dehn twists $\tau_{\b_1}$ and $\tau_{\b_2}$
about disjoint curves $\b_1$ and $\b_2$ on a surface $F$, the
triangles determined by the points $x$, $\tau_{\b_1}^n\cdot x$,
$\tau_{\b_2}^{-n}\cdot x$ contradict Gromov's condition: the legs
of
this triangle are given by the \tec geodesics whose corresponding
\tec maps from $x$ are described explicitly in \cite{MM}, and the
distances between points on the legs are estimated from below in
terms of estimates of relevant extremal lengths.

We organize our discussion as follows. In \S2, we recall the
background information we will need, and set the notation. In \S3
we state and prove our main result.

The authors are indebted to M\.~Kapovich for bringing this
question
to their attention.

\newpage
\sub{\S2. Background and Notation}

\mpage
\sub{\S2.1. \tec space, metric, maps} Let $M$ be a closed
$C^\infty$ surface of genus $g\ge2$; everything in this note
extends to
punctured surfaces with only additional notation, so we
concentrate
on the closed surface case. We consider the \tec space $T_g$ with
the \tec metric $d(\cd,\cd)$. Recall that points in \tec space
are
equivalence classes of Riemann surface structures $S$ on $M$, the
structure $S_1$ is equivalent to the structure $S_2$ if
there is a homeomorphism $h:M\to M$, homotopic to the identity,
which is a conformal map of the structures $S_1$ and $S_2$.

We define the \tec distance $d(\{S_1\},\{S_2\})$ by
$$
d(\{S_1\}, \{S_2\}) = \f12\log\inf_h K(h)
$$
where $h: S_1\to S_2$ is a quasiconformal homeomorphism homotopic
to the identity on $M$ and $K[h]$ is the maximal dilatation of
$h$.
This metric is well-defined, so we may unambiguously write
$S_1$ for $\{S_1\}$.

An extraordinary fact about this metric is that the extremal
maps,
known as \tec maps, admit an explicit description, as does the
family of maps which describe a geodesic.

Specifically, let $q\in\qd(S)$ denote a holomorphic quadratic
differential on $S$. A horizontal (resp\. vertical) trajectory is
an arc along which $q(z)dz^2>0$ (resp\. $q(z)dz^2<0$) except at
the
zeros of $q$. A trajectory is critical if it passes through a
critical point; otherwise it is regular. If $z$ is a local
parameter near $p\in S$ with $q(p)\neq0$ and $z(p)=z_0$, then
$w=\int^z_{z_0}q(z)^{1/2}dz$ is the natural parameter $q$ near
$p$.
The line element $|q(z)|^{1/2}|dz|$ defines the $q$-metric on
$S$.

\tecs theorem asserts that if $S_1$ and $S_2$ are distinct
points in $T_g$, then there is a unique quasiconformal
$h : S_1\to S_2$ with $h$ homotopic to the identity on $M$ which
minimizes
the maximal dilatation of all such $h$. The complex dilatation of
$h$ may be written $\mu(h)=k\f{\bar q}{|q|}$ for some non-trivial
$q\in\qd(S_1)$ and some $k$, $0<k<1$, and then
$$
d(S_1,S_2) = \f12\log(1+k)/(1-k).
$$
Conversely, for each $-1<k<1$ and non-zero $q\in\qd(S_1)$, the
quasiconformal homeomorphism $h_k$ of $S_1$ onto $h_k(S_1)$,
which
has complex dilatation $k\bar q/|q|$, is extremal in its homotopy
class. Each extremal $h_k$ induces a quadratic differential
$q'_k$
on $h_k(S_1)$, with critical points of $q$ and $q'_k$
corresponding
under $h_k$; furthermore, to the natural parameter $w$ for $q$
near $p\in S_1$ there is a natural parameter $w'_k$ near $h_k(p)$
so that
$$
\Ree w'_k = K^{1/2}\Ree w\quad\text{and}\quad
\im w'_k = K^{-1/2}\im w,
$$
where $K=(1+\e)/(1-k)$.

The map $h_k$ is called the \tec extremal map determined by $q$
and
$k$;the differential  $q$ is called the initial differential and
the differential $q_k$ is called the
terminal differential.   We can assume all quadratic
differentials are normalized in the sense that $$||q||=\int
|q|=1.$$ The \tec geodesic segment between $S_1$ and $S_2$
consists of all points $h_s(S_1)$ where the $h_s$ are
\tec maps on $S_1$ determined by the quadratic differential
$q\in\qd(S_1)$ corresponding to the \tec map $h : S_1\to S_2$
and $s\in[0, \|\mu(h)\|_\infty]$.

The mapping class group $Diff^+(M)/Diff_0(M)$ acts on $T_g$.  If
$\{U_\alpha,z_\alpha\}$ is an atlas defining the Riemann surface
structure $S$, and $f$ is a diffeomorphism of $M$, then $f\cdot
S$ is the Riemann surface structure defined by the atlas
$\{f(U_\alpha), z_\alpha\circ f^{-1}\}$.  The map $f:S\to f\cdot
S$ is then a conformal map between these two structures.

\bpage
\sub{\S2.2. Modulus, Extremal length, Jenkins-Strebel
Differentials, Dehn twists}\newline
The modulus of a flat cylinder $C$ of
circumference $l$ and height $h$ is $\bmod(C)=h/l$. For a simple
closed curve $\g\subset M$, we define the modulus $\bmod_S(\g)$
of
$\g$ to be the supremum of the moduli of all cylinders embedded
in
$M$ with core curve isotopic to $\g$.

The extremal length $\ext_S(\g)$ of a curve $\g$ on a surface $M$
is defined to be $$\sup_\rho(l_\rho([\g]))^2/A_\rho,$$ where
$\rho$ ranges over all conformal metrics on $S$ with area
$0<A_\rho<\infty$ and $l_\rho([\g])$ denotes the infimum of
lengths
of simple closed curves homotopic to $\g$. One shows that
$\ext_S(\g)=1/\bmod_S(\g)$.

Kerckhoff \cite{K} has given a characterization of the \tec
metric
$d(S_1, S_2)$ in terms of the extremal lengths of
corresponding curves on the surfaces. He proves
$$
d(S_1, S_2) = \f12\log\sup\lm_\g\
\f{\ext_{S_1}(\g)}{\ext_{S_2}(\g)}\tag2.1
$$
where the supremum ranges over all simple closed curves on $M$.

Jenkins \cite{J} and Strebel \cite{Str} proved the existence of
quadratic differentials $q\in\qd(S)$ with some prescribed
trajectory topology. Specifically, they (see \cite{Str}, e.g.)
showed that one could specify $m$ disjoint simple loops $\vecs\g
m$,
with $1\le m\le 3g-3$, on $S$ representing distinct non-trivial
free homotopy classes, and $m$ positive numbers $\vecs Mm$, and
that then one could find a unique (up to scalar multiple)
quadratic
differential $Q=Q(z)dz^2\in\qd(S)$ with the following property:
if
$S'$ is the result of removing the critical trajectories of
$Q(z)dz^2$ from $S$, then $S'$ is the union of annuli $\vecs Am$
with $A_j$ homotopically equivalent to $\g_j$ and the modulus of
$A_j$ was $M_j$, up to some fixed (independent of $j$) scalar
multiple. Further $S-S'$ is the union of a finite number of
analytic arcs, the smooth pieces of the critical trajectories.

Consider a point  $S\in T_g$ and consider the effect of a
Dehn twist $\tau_\a$ about a curve $\a\subset M$ yielding a point
$\tau_\a \cdot S\in T_g$. It is natural to
ask for a
characterization of the \tec map $h:S\to \tau_\a\cdot S$, or
more generally, for a characterization of the \tec map
$h_n$ from $S\to\tau_\a^n\cdot S$ in terms of the data
$\a$, $S$
and $n\in\BZ$. This was described by Masur and Marden \cite{MM}
as
follows. Let $q_\a=q_\a(z)dz^2$ denote the Jenkins-Strebel
differential determined, as above, by $\a\in M$, and
suppose that $\a\subset S$ has modulus $R$. Set
$$
\align
M  &= (\log R)/2\pi\\
\s_n  &=\tan^{-1}(2M/n)\\
\intertext{and}
k_n  &= \f{|n|/2M}{\(1+\(\f n{2M}\)^2\)^{1/2}}.
\endalign
$$
Then \cite{MM} the extremal map $h_n : S\to \tau_\a^n\cdot S$ is
the \tec map determined by $[\exp(-i(\s_n+\pi))]q_\a$ and
$k_n$. Furthermore, if we pull back the terminal quadratic
differential $q_\alpha'$ on
$\tau_\a^n\cdot S$ to $S$ via the conformal map $\tau_\a^n$, then
the
pull-back differential $(\tau_\a^n)^*q_\alpha'$ satisfies
$$(\tau_\a^n)^*q_\alpha'=e^{i\th}q_\alpha\tag 2.2$$ in particular
the metrics
$|q_\alpha|$ and
$|(\tau_\a^n)^*q_\alpha'|$ agree.

\bpage
\sub{\S2.3. Gromov hyperbolicity} Let $X$ be a geodesic metric
space, that is, a metric space $(X,d)$ where every pair of points
$x$, $y\in X$ can be connected by the isometric image of the
segment $[0, d(x,y)]$. In such a space, we can define the notion
of
a triangle with vertices $x$, $y$ and $z\in X$ to be the union of
geodesic segments $[xy]$, $[yz]$, and $[xz]$ connecting $x$ and
$y$, $y$ and $z$, and $x$ and $z$, respectively. Naturally, \tec
space with the \tec metric is a geodesic metric space.

Gromov (see \cite{GdlH}) introduced a notion of when such a space
would share a number of qualitative properties with hyperbolic
space, his definition now being commonly called ``Gromov
hyperbolicity''. We will say that

\bpage
\sub{Definition 2.1} The geodesic metric space $X$ is Gromov
hyperbolic if
$$
\align
&\text{There is a number $\d\ge0$ so that for every triangle
$\Delta=[xy]\cup[yz]\cup[xz]$ and}\tag"($*$)"\\
&\text{every $u\in[xy]$, we have } d(u, [yz]\cup[zx])\le\d.
\endalign
$$

Hyperbolic space, (Riemannian) negatively curved manifolds,
trees,
Euclidean strips, free groups with the word metric and spheres
are
easily shown to be Gromov hyperbolic. On the other hand, the
fundamental group of a non-compact finite volume hyperbolic
$n$-manifold with $n\ge3$, equipped with the word metric, is not
hyperbolic, because of the large rank (parabolic) abelian
subgroup
stabilizing a point at infinity (cusp).

\newpage
\sub{\S3. Main Theorem}

\mpage
The goal of this section is to prove

\bpage
\proclaim{Theorem 3.1} \tec space with the \tec metric is not
Gromov hyperbolic.
\endproclaim

\bpage
\sub{Proof} We consider a sequence of triangles $T_n$ so that
there
does not exist a $\d\ge0$ with condition ($*$) (in
Definition~2.1)
holding for all $T_n$.

All the triangles $T_n$ will have a common vertex $x_0\in T_g$,
chosen arbitrarily. The other vertices of the triangle $T_n$ are
the points $y_1=\tau_{\b_1}^n\cdot x_0$ and
$y_2=\tau_{\b_2}^{-n}\cdot x_0$,
where $\b_1$ and $\b_2$ are disjoint simple closed curves on the
surface $M$ of genus $g>1$.

We wish to estimate the \tec distance from a point $y\in[y_1y_2]$
to the other legs $[x_0y_1]$ and $[x_0y_2]$. To this end, we let
$J_1dz^2\in\qd(x_0)$ be the Jenkins-Strebel differential with
core
curves homotopic to $\b_1$, and we suppose that the union of its
regular trajectories determine an annulus of modulus $R_1$. We
let
$M_1=(\log R_1)/2\pi$, $\tan\tau_1=2M/n$, and
$k_1=|n|(2M_1)^{-1}(1+(n/2M_1)^2)^{-1/2}$, so that the \tec map
from
$x_0$ to $y_1$ is determined by $\exp(-i(\tau_1+\pi))J_1$ and
$k_1$.

Let $\g_1$ be a simple closed curve on $M$ which crosses $\b_1$
but
not $\b_2$, and let $\g_2$ be a simple closed curve on $F$ which
crosses $\b_2$ but not $\b_1$. Then we claim

\bpage
\proclaim{Lemma 3.2} For $x\in[x_0y_1]\subset T_n\subset T_g$,
the extremal length, $\ext_x(\g_2)$, of $\g_2$ on $x$ is bounded
independently of $n$.
\endproclaim

\bpage
\sub{Proof of Lemma 3.2} We begin with some more notation.
Consider
a quadratic differential $q\in\qd(x_0)$ and the associated
singular
flat Euclidean metric $|q|$. For a $|q|$-geodesic segment $\a$,
let
the horizontal and vertical $q$-lengths of $\a$ be denoted
$$
\split
h_q(\a)  &= \int\lm_\a |\Ree q^{1/2}|\\
v_q(\a)  &= \int\lm_\a |\im q^{1/2}|.
\endsplit
$$
Then
$$
|\a|_q = (h_q(\a)^2 + v_q(\a)^2)^{1/2},\tag3.1
$$
where $|\a|_q$ is the $q$-length of $\a$. We observe that under
the
\tec map determined by $q$ and $K$ with terminal quadratic
differential $q'$, we'll have the arc $\a$ remaining a
$q'$-geodesic arc and
$$
h_{q'}(\a) = K^{1/2}h_q(\a),
v_{q'}(\a) = K^{-1/2}v_q(\a)\ \text{ and }\
|\a|^2_{q'} = Kh_q(\a)^2+K^{-1}v_q(\a)^2.\tag3.2
$$

Of course, for fixed $h_q(\a)$ and $v_q(\a)$, equation (3.2)
expresses $|\a|_{q'}$ as a convex function of $K>0$.

We now specialize to the case in the statement of the lemma,
where
$J_1\in\qd(x_0)$ determines the \tec geodesic arc
$[x_0y_1]\subset T_g$ and $J_1'$ is the terminal differential on
$y_1$.  Since $\tau_\a^n(\b_1)=\b_1$, (2.2) implies
$$|\b_1|_{J_1}=|\b_1|_{J_1'}.$$ The convexity of $|\b_1|$
in
$K$ along $[x_0y_1]$ forces
$|\b_1|_{J_x}<|\b_1|_{J_1}=|\b_1|_{J_1'}$ for any of the
quadratic
differentials $J_x\in\qd(x)$ associated to the \tec geodesic
segment
$[x_0y_1]$ and any $x\in[x_0y_1]^0$. On the other hand, because a
\tec map is area preserving, this forces
$$
\bmod_x(\b_1) > \bmod_{x_0}(\b_1) = \bmod_{y_1}(\b_1)\tag3.3
$$
where $\bmod_x(\b_1)$ refers to the modulus of the $\b_1$ annulus
on $x\in[x_0y_1]$.

We use (3.3) in considering an alternative description of the
\tec
map between $x_0$ and $x\in[x_0y_1]$. Specifically, by the same
technique of proof as that for Lemma~2.1 in \cite{MM} (see also
the
statement for the annulus in \cite{MM; \S1.3}), we can represent
the \tec map between $x_0$ and $x\in[x_0y_1]$ as $T_\a\circ S_a$
where $T_\a$ is a ``partial'' Dehn twist of the initial
Jenkins-Strebel annulus by an angle $2\pi\a$ and $S_a$ is a
radial
expansion or (possibly) contraction of that annulus: we observe
however that by (3.3), the map $S_a$ is {\it always\/} an
expansion.

Thus, we can build a model of any terminal Jenkins-Strebel
differential $J_x\in\qd(x)$ with $x\in[x_0y_1]$ as given by an
operation of conformal plumbing followed by a partial Dehn twist,
as follows. We cut the conformal cylinder along a core curve. We
then glue in one cylinder to each edge of the cut, again leaving
a
pair of boundary components. Finally, we glue these free edges
together after twisting by some angle.

The homotopy class of $\g_2$ is represented by a union of
geodesic
segments on the boundary of the Jenkins-Strebel annulus for
$J_1$.
Therefore, we can find an annulus $A_2$, embedded around $\g_2$,
and also disjoint from the core curve along which our initial cut
(of the previous paragraph) is made. That annulus $A_2$ will be
unaffected by the plumbing and twisting, and so we can conclude
that
for all $x\in[x_0y_1]$ for which $x=T_a\circ S_\a x_0$, we can
find
an embedded annulus $A_2$ about $\g_2$ of modulus bounded
uniformly
away from zero, independently of $n$.

Thus the extremal length of $\g_2$ is then uniformly bounded
above,
independently of $n$, concluding the proof of the lemma. \qed

\bpage
\sub{Remark} The lemma of course holds with $\g_1$ and $[x_0y_2]$
in place of $\g_2$ and $[x_0y_1]$, by an interchange of notation
in the proof.

\bpage
\sub{Conclusion of the proof of Theorem 3.1} Now consider the
\tec
geodesic arc $[y_1y_2]$. The \tec map from $y_1$ to $y_2$ is
given
by taking a negative twist $n$ times about $\b_1$ and about
$\b_2$.  Consider the Strebel differential $Q\in\qd(y_1)$ of two
annuli with core curves homotopic to $\b_1$ and $\b_2$, of equal
moduli $R$ (see \cite{Str}).  Let $M$, $\s_n$ and $k_n$ be as
in
\S2.2; then the \tec map from $y_1$ to $y_2$ is determined by
$\exp(-i(\s_n+\pi))Q$ and $k_n$.  Let $Q'$ be the terminal
differential on $y_2$.

By Lemma 3.2 and the fact that $Q$ is a competing metric in the
definition of extremal length, we have
$$
|v_Q(\g_2)|\leq |\gamma_2|_Q\leq
\ext_{y_1}(\g_2)^{1/2}=O(1).\tag3.4
$$
Since $y_2=\tau_{\beta_2}^{-n}\cdot x_0$, we have
$$\ext_{y_2}(\gamma_2)=\ext_{x_0}(\tau_{\beta_2}^n(\gamma_2)).$$
Since $\tau_{\beta_2}(\gamma_2)$ crosses $\beta_2$ $n$ times,
there
is a constant $c_0>0$ so that $$\ext_{y_2}(\g_2)\ge c_0n^2.$$
Moreover, since
we can always compare any two normalized metrics on the fixed
 surface $y_2$, conformally equivalent to $x_0$, we find that
$$
|\g_2|_{Q'}\ge cn\tag3.5
$$
for some $c>0$.

Next, since
$$
k_n = \(1 + \(\f{\log R}{\pi|n|}\)^2\)^{-1/2}
$$
we see that
$$
K_n = \f{1+k_n}{1-k_n} \asymp n^2 \tag3.6
$$ where $a\asymp b$ if their ratio is bounded above and below
away from $0$.
Then, applying (3.4), (3.5) and (3.6) to the identity
$$
K_nh_Q(\g_2)^2 + K_n^{-1}v_Q(\g_2)^2 = |\g_2|_{Q'}^2
$$
yields
$$
h_Q(\g_2) > c_2 > 0.\tag3.7
$$

Next, we observe that $-Q$ is the terminal quadratic differential
on $y_1$ for the \tec map from $y_2$ to $y_1$, with initial
differential $-Q'$. Then the same argument as above shows that
$h_{-Q'}(\g_1)>c'_2>0$, independently of $n$. We can then apply
formula (3.2) again to conclude that $h_{-Q}(\g_1)>c_3n$ for some
$c_3>0$, which, of course, is equivalent to
$$
v_Q(\g_1) > c_3n.\tag3.8
$$

Finally, consider the point $y_*\in[y_1y_2]$ determined by the
\tec
map defined by $Q$ with $K^{1/2}=\sqrt n$; let $Q_*\in\qd(y_*)$
denote the terminal differential. Then (3.7) and (3.8), along
with
the relationship (3.2) show that
$$
\split
|\g_2|_{Q_*} &\ge h_{Q_*}(\g_2)\ge c_2\sqrt n\qquad\text{and}\\
|\g_1|_{Q_*} &\ge v_{Q_*}(\g_1)\ge c_3\sqrt n.
\endsplit
$$
Since $Q_*$ is a competing metric for extremal length,
$\ext_{y_*}(\g_i)\ge |\g_i|_{Q_*}^2>c_4n$.

Finally, we apply Kerckhoff's formula (2.1) and Lemma 3.2 to
estimate
the \tec distance $d([x_0y_1],y_*)$: we see that since
$\ext_{y_*}(\g_2)>c_4n$ while $\ext_x(\g_2)<c_5$ for
$x\in[x_0y_1]$, then (2.1) forces
$d(x,y_*)>\f12\log(c_5^{-1}c_4n)$. Since an analogous
estimate
holds for $d([x_0y_2],y_*)$, we see that the defining condition
$(*)$ of Definition~2.1 of Gromov hyperbolicity does not hold.

\end